\documentstyle{amsppt}
\vsize= 53pc
\hsize= 35pc
\NoBlackBoxes

\topmatter 
\title Detecting quasiconvexity: algorithmic aspects \endtitle
\rightheadtext {Detecting quasiconvexity}
\author Ilya Kapovich \endauthor
\address Department of Mathematics, City College of the City University of New York, Convent Avenue at 138-th Street, New York, NY10031\endaddress
\email ilya\@groups.sci.ccny.cuny.edu\endemail
\abstract 
 The main result of this paper states that for any group $G$ with an automatic structure $L$ with unique representatives one can construct a uniform partial algorithm which detects L-rational subgroups and gives their preimages in $L$. This provides a practical, not just theoretical, procedure for solving the occurrence problem for such subgroups.	 
\endabstract
\endtopmatter
%\mag\magstep1
%\hsize=6.5 true in

\document

\head 1. Generalized word problem and rational structures on groups\endhead
The goal of this paper is to highlight connections between the theory of automatic groups and the generalized word problem and to demonstrate certain additional advantages of the class of automatic groups over the class of combable groups. We assume that the reader is familiar with the theory of automatic groups, regular languages and combable groups.  Although some of the important definitions will be given, the reader is referred to [ECHLPT] for further details. A good overview of the theory of automatic groups can also be found in [BGSS]. We take for granted some basic facts about word hyperbolic groups and their connections with the theory of automatic groups. Here our main references are [Gr], [ABCFLMSS], [ECHLPT], [BGSS] and [GS]. An important discussion about combable groups can also be found in [A], [AB] and [N]. 
The author is grateful to the referee for greatly simplifying the proof of Proposition 1 and to Gilbert Baumslag for his help in writing this paper.

Recall that if $G$ is a recursively presented finitely generated group given by a presentation $$G=<x_1,..,x_n | r_1,..,r_m,..>  \eqno (1)$$
and $H$ is a subgroup of $G$ then we say that $G$ has {\it solvable generalized word problem} with respect to $H$ if there is an algorithm which, for any word $w$ in the generators $x_1,..,x_n$, decides whether or not $w$ represents an element of $H$.

Equivalently, $G$ has solvable generalized word problem  with respect to $H$ if the set $\phi^{-1}(H)$ is a recursive subset of the free group $F(x_1,..,x_n)$ where $\phi\colon F(x_1,..,x_n)\rightarrow G$ is the natural epimorphism associated with the presentation (1).
It is not hard to see that this definition does not depend on the choice of a presentation of $G$ with a finite number of generators.
\smallskip
If $G$ has solvable generalized word problem with respect to the trivial subgroup $H=1$ then $G$ is said to have {\it solvable word problem }.
\smallskip
The concept of generalized word problem goes back to the work of W.Magnus [M] where he proved that if $M$ is a subgroup of a one-relator group $G=<x_1,..,x_n | R=1>$, generated by any subset of the generating set, then $G$ has solvable generalized word problem with respect to $M$. 
We should stress that the mere knowledge that $G$ has solvable generalized word problem with respect to $H$ does not yet give one an effective procedure for determining whether or not a particular word in the generators of $G$ represents an element of $H$. That is, for practical purposes it is not enough to know that the required algorithm exists, it is necessary to be able to find it.

\definition {Definition} Let $G$ be a finitely generated group. A {\it rational structure} on $G$ is a pair $(L,A)$, where $A$ is a finite generating set for $G$, closed under taking inverses, and $L$ is a regular language over the alphabet $A$ with $\pi(L)=G$.  
(Here $\pi\colon A^{\ast}\rightarrow G$ is the natural monoid homomorphism from the free monoid $A^{\ast}$ on $A$ to the group $G$.)
We say that $(L,A)$ is a rational structure {\it with uniqueness} if for any element $g\in G$ there is a unique element $w\in L$ such that $\pi(w)=g$.
\enddefinition

To distinguish the monoid multiplication in $A^{\ast}$ from the group multiplication in $G$ we will use $w_1w_2$ for the former and $g_1\cdot g_2$ for the latter. We will also, when no confusion is possible, denote the image $\pi(w)$ in $G$ of a word $w\in A^{\ast}$ by $\overline w$.

The notion of a rational structure on a group was first closely investigated by R.Gilman in [Gi4] who used the term "rational cross-sections" for rational structures with uniqueness. A lot of important facts about groups with rational structures, which later were incorporated into the theory of automatic groups, can be found in the early works of R.Gilman [Gi1],[Gi2],[Gi3] and [Gi4].

One should think of a rational structure $(L,A)$ for $G$ as of a regular collection of normal forms for the elements of $G$. This point of view is propagated in particular in the work of H.Short [Sho] where he investigates different types of rational structures on finitely generated groups.

\proclaim {Lemma 1} Suppose $G=<x_1,..,x_n | r_1,..,r_m,..>$ is a recursively  presented group which admits a rational structure with uniqueness $(L,A)$ where $A=\lbrace x_1^{\pm 1},..,x_n^{\pm 1}\rbrace$.
Then $G$ has solvable word problem.
\endproclaim
\demo {Proof}
There is a unique word $w_{\ast}\in L$ representing the identity element of $G$. 

Since $G$ is recursively presented, the normal closure $N$ of the set $\lbrace r_1,..,r_m,..\rbrace$ in the free group $F=F(x_1,..,x_n)$ is recursively enumerable.
Suppose now $v$ is a freely reduced word over $x_1,..,x_n$.
The language $L$ is regular and so it can be recursively
enumerated. Thus one can also recursively enumerate the set $N\times L$.
We now inspect each pair $(n,w)$ in $N\times L$ and check if $n=v^{-1}\cdot w$ in $F$. There is a unique $w\in L$ such that $\overline w=\overline v$ in $G$ that is $v^{-1}\cdot w\in N$.
When we find a pair $(n,w)$ with $n=v^{-1}\cdot w$ in $F$, we conclude that $\overline v=1$ in $G$ if $w=w_{\ast}$ and that $\overline v\not=1$ otherwise. This completes the proof.
\enddemo

\proclaim {Corollary 2} (see [ECHLPT]) If $G$ is an asynchronously automatic group (in the sense of [ECHLPT]) then $G$ has solvable word problem.
\endproclaim
\demo {Proof} Indeed, by Theorem 7.3.2 and Theorem 7.3.4 of [ECHLPT], the group $G$ is finitely presentable and for any finite generating set it possesses an asynchronous automatic structure with uniqueness. Thus the statement follows from Lemma 1.
\enddemo

The following definition, due to S.Gersten and H.Short [GS], is of importance here.
\definition {Definition} Let $(L,A)$ be a rational structure on $G$ and let $H$ be a subgroup of $G$.
\roster
\item "(a)" The subgroup $H$ is said to be $L$-{\it rational} if its full preimage in $L$, $L_H=L\cap \pi^{-1}(H)$, is a regular language.
\item "(b)" The subgroup $H$ is said to be $L$-{\it quasiconvex} if there is a constant $K>0$ such that for any $w\in L$ with $\overline w\in H$ and for any initial segment $w_t$ of $w$ there is a word $u_t\in A^{\ast}$ of length at most $K$ such that $\overline {w_tu_t}\in H$.
\endroster
\enddefinition

An important observation made by S.Gersten and H.Short in [GS] asserts that, adopting the notation above, $H$ is $L$-quasiconvex if and only if $H$ is $L$-rational. It is shown in the same paper that an $L$-quasiconvex subgroup $H$ is always finitely generated.

\proclaim {Lemma 2} Let $(L,A)$ be a rational structure with uniqueness for a finitely generated recursively presentable group $G$ and let $H$ be an $L$-rational subgroup of $G$.
Then the generalized word problem for $G$ with respect to $H$ is solvable.
\endproclaim
\demo {Proof}
Indeed, if $v\in A^{\ast}$ then as the proof of Lemma 1 shows, one can algorithmically find the normal form $w\in L$ of $\overline v$.
Since $\pi^{-1}(H)\cap L=L_H$ is regular and is accepted by some finite state automaton $M_H$, it remains only to check whether $M_H$ accepts $w$ or not.
\enddemo

One may think that Lemma 2 completes the discussion about the generalized word problem and rational structures on groups.
However, as we pointed out earlier, the question of effectively finding the desired algorithm requires some further investigation.
We will concentrate on the following questions.
\smallskip
{\bf Question 1} Suppose a group $G=<x_1,..,x_n|r_1,...,r_m>$ is given by a finite presentation. Suppose further that a rational structure with uniqueness $(L,A)$, where $A=\{x_1^{\pm 1},..,x_n^{\pm 1}\}$, is given explicitly by a finite state automaton $M$ recognizing $L$. 

Is there a uniform (on $H$) procedure for finding the language $L_H=\pi^{-1}(H)\cap L$ if the subgroup $H$ of $G$, given by a finite set of generators, is known to be $L$-rational?
\smallskip
{\bf Question 2} Suppose $G=<x_1,..,x_n|r_1,...,r_m>$ is given by a finite presentation and that $H=gp(\overline {v_1},...,\overline {v_t})$ is given by a generating set $v_1,..,v_t$. Suppose that we know that there exists a rational structure with uniqueness for $G$ such that $H$ is rational.

Is there a uniform (on $H$) procedure for finding the automata $M$ and $M_H$ such that $(L=L(M),A)$ is a rational structure with uniqueness for $G$ and $L(M_H)=L_H=\pi^{-1}(H)\cap L$?
\smallskip
We will see that the answers to both these questions are positive if we restrict ourselves to the class of automatic structures.

\head 2. Detecting rational subgroups of automatic groups\endhead

\proclaim {Proposition 1} Suppose $G=<x_1,..,x_n | r_1,..,r_m>$ is a group given by a finite presentation. Suppose $(L=L(M),A=\{x_1^{\pm 1},..,x_n^{\pm 1}\}, M_{=},M_{x_1},...,M_{{x_n}^{-1}})$ is an automatic structure with uniqueness for $G$. 

Then there is a uniform (on $H$) partial algorithm which, for a subgroup $H=gp(\overline {v_1},..,\overline {v_t})$ given by a finite generating set ${v_1},..,{v_t}$, will
\roster
\item "(a)" eventually stop and produce the automaton $M_H$, recognizing 

$L_H=\pi^{-1}(H)\cap L$, when $H$ is $L$-rational
 
\item "(b)" run forever when $H$ is not $L$-rational.
\endroster
\endproclaim

\demo {Proof}

Recall that a {\it Schreier diagram} of $H$ in $G$ with respect to the generating set $A=\{x_1^{\pm 1},..,x_n^{\pm 1}\}$ is the labeled oriented graph $\Gamma (G,H,A)$ whose vertices are cosets $Hg$, $g\in G$ and which has an oriented edge
$(Hg, Hga)$ labeled by $a$ for each vertex $Hg$ and every $a\in A$.
Thus a word in the generators $A$ represents an element of $H$ if and only if it is a label of a cycle starting at the vertex $H=H\cdot 1$ of $\Gamma (G,H,A)$. 
It follows from the definition of $L$-quasiconvexity that the subgroup $H$ is $L$-quasiconvex if and only if there is
$K>0$ such that any word in $L$ representing an element of $H$
is a label of a cycle in  $\Gamma (G,H,A)$ starting from $H$ and contained  in the ball of radius $K$ around $H$ in  $\Gamma (G,H,A)$.
\smallskip
We will now describe the required algorithm.
We may assume that the unit element $1\in G$ is represented by the empty word $\varepsilon\in L$. If not, we find the word $w_{\ast}\in L$ such that $\overline w=1$. Then $L'=(L-{w_{\ast}})\cup {\epsilon}$ is an automatic language with uniqueness for $G$. We construct the automaton recognizing $L'$ and modify the comparison automata $M_=, M_{x_1},..,M_{x_n^{-1}}$ accordingly.
\smallskip
{\bf Step 0} First, for each word $v_i$, $i=1,..,n$ construct the automaton $M_{v_i}$ which accepts all pairs $(w,u)\in L\times L$ such that $\overline {wv_i}=\overline u$.
Do the same for the inverse of $v_i$.
\smallskip
{\bf Step 1} Use the Todd-Coxeter method (see ch. 5 of [Si]) to enumerate right cosets of $H$ in $G$. The method produces a sequence of finite labeled graphs $X_i$ with labels from $A$. These $X_i$ approximate bigger and bigger balls in the Schrier diagram $\Gamma (G,H,A)$. Each $X_i$ has a distinguished basepoint standing for $H$ and the label of any cycle in $X_i$, starting at the basepoint, represents an element of $H$.

\smallskip
Each $X_i$ can be turned into a finite state automaton by taking the basepoint to be its initial and terminal state. Thus a word over $A$ is accepted by $X_i$ if and only if this word is the label of a cycle in $X_i$ starting at the basepoint.
Therefore any word accepted by $X_i$ represents an element of $H$.
Moreover, for any integer $k>0$ the balls of radius $k$ in the graphs $X_i$ eventually stabilize (as $i$ tends to infinity) and map isomorphically (preserving the labels and the basepoints) to the ball of radius $k$ in $\Gamma (G,H,A)$.
Therefore if $H$ is $L$-quasiconvex, then for $i$ large enough every word from $L$ representing an element of $H$ will be accepted by $X_i$.
\smallskip
{\bf Step 2} Compute the automaton accepting $L_i$, the intersection of $L$ and the language accepted by $X_i$.
For each generator $v_j$ of $H$ check if $L_i$ is "stable" under right multiplication by $v_j$.
That is, compute the intersection $L_i\times L_i\cap L(M_{v_j})$ and its projection $N_i$ on the first coordinate.
Then check if $N_i=L_i$.
Do the same for $v_j^{-1}$.
If the answers are yes for each $j=1,..,t$ then stop the procedure. The output of the algorithm is the language $N=L_i$.
If the answer is no, go to Step 1 and increase $i$ by 1.
\smallskip

We claim that if the algorithm stops and produces the language $N$ then $H$ is $L$-rational and $N=L_H$.
Indeed, by construction every word from $N$ lies in $L$ and represents an element of $H$. Also, the empty word $\varepsilon$ belongs to $N$. It follows now from the description of Step 2 that for any $k=1,..t$ $$1\in \overline N, \qquad  \overline {Nv_k}\subset \overline N \quad\hbox{and}\quad \overline {Nv_k^{-1}}\subset \overline N.$$ Thus $\overline N= H$ and, since $L$ is an automatic language with uniqueness, $N=L_H$.
\smallskip
On the other hand, if $H$ is $L$-rational, it is $L$-quasiconvex and, as it was observed in the description of Step 1, for some $i$ every word from $L$ representing an element of $H$ will be accepted by $X_i$. It is clear that for such $i$ we have $L_i=L_H$ and the language $L_i$ is stable under right multiplication by $v_j$ and $v_j^{-1}$, $j=1,..,t$.
Thus the algorithm stops and produces the language $N=L_H$.
\smallskip
This completes the proof of Proposition 1.
 
\enddemo

Proposition 1 implies using the notations above, that the set of $L$-rational subgroups of $G$ is recursively enumerable.

\proclaim {Proposition 2} Suppose $G=<x_1,..,x_n | r_1,..,r_m>$ is a group given by its finite presentation. 
There is a uniform partial algorithm which, for a subgroup $H=gp(\overline {v_1},..,\overline {v_t})$, given by a finite generating set $v_1,..,v_t$, will

 \roster 
\item "(a)" eventually stop and produce an automatic structure with uniqueness $$(L=L(M),A,M_{=},M_{x_1},..,M_{{x_n}^{-1}})$$ for $G$ and a finite state automaton $M_H$ such that $L(M_H)=\pi^{-1}(H)\cap L$ when $H$ is rational for some automatic structure on $G$;
\item "(b)" run forever if there is no automatic structure for $G$ such that $H$ is rational with respect to it.
\endroster
\endproclaim

\demo {Proof} 

The statement easily follows from Proposition 1 and Theorem 5.2.4 of [ECHLPT] which provides a uniform partial algorithm for detecting automatic groups and building automatic structures on them.
\enddemo

Most results concerning automatic groups have counterparts for combable groups. The only exceptions known to us  are the theorem of R.Gilman ([Gi4]) which asserts that infinite groups with rational structures with uniqueness have elements of infinite order, the theorem of Epstein and Holt [ECHLPT] which states that nilpotent groups are not asynchronously automatic and Theorem 5.2.4 of [ECHLPT] which provides a uniform partial algorithm for detecting automatic groups.
Proposition 1 and Proposition 2, which also have no analogs for asynchronously automatic and combable groups, seem to be interesting additions to this list. The author is grateful to the referee for greatly simplifying the proof of Proposition 1.

\head  3. Hyperbolic Groups: an alternative approach\endhead

 Proposition 1 provides a general algorithm for detecting quasiconvex subgroups of automatic groups.
There is an alternative approach for detecting quasiconvex subgroups of word hyperbolic groups. This approach was suggested to the author by P.Papasoglu and it is indeed very much in the spirit of Papasoglu's paper [P].

Recall that word hyperbolic groups are $ShortLex$-automatic for any finite generating set (see [ECHLPT], Theorem 3.4.5).

A metric space $(X,d)$ is termed {\it geodesic} if any two points can be joined by a path whose length is equal to the distance between these points. Such a path, when parametrized by arc-length, is also called a {\it geodesic}. A naturally parametrized path is called $\lambda$-{\it quasigeodesic} if for any points $x$ and $y$ on it 
 $$s\le \lambda d(x,y)+\lambda$$ where $s$ is the length of the segment of this path between $x$ and $y$.  A naturally parametrized path is called $\mu$-{\it local} $\lambda$-{\it quasigeodesic} if any segment of this path of length at most $\mu$ is $\lambda$-quasigeodesic.
If $G$ is a group and $A$ is a finite generating set of $G$, we denote the Cayley graph of $G$ with respect to $A$ by $\Gamma(G,A)$. The word metric on $\Gamma (G,A)$ is denoted $d_A$. 

We collect some useful facts about word hyperbolic groups in the following statement.

\proclaim {Proposition 3} Let $G$ be a word hyperbolic group.
 
(i)  Let $H$ be a subgroup of $G$. Then the following conditions are equivalent.

\roster
\item "(1)" For some (and so for any) finite generating set $A$ of $G$ and some (and so for any) automatic structure $(L,A)$ on $G$ the subgroup $H$ is $L$-rational.
\item "(2)" $H$ is finitely generated and for some (and so for any) finite generating set $B$ of $H$ and for some (and so for any) finite generating set $A$ of $G$ the inclusion $(H,d_B)\subset (G,d_A)$ is a quasiisometry, that is there is some $C>0$ such that

$$(1/C)d_B(h,1)-C\le d_A(h,1)\le Cd_B(h,1)+C$$ for each $h\in H$.

We call $C$ a {\it distortion constant} for $H$.

\item "(3)" For some (and so for any) finite generating set $A$ of $G$ the subgroup $H$ is a quasiconvex subset of $(\Gamma (G,A),d_A)$, that is there is an $\varepsilon >0$ (called the {\it quasiconvexity constant} for $H$) such that for each $h\in H$ and for any $d_A$-geodesic path $w$ from $1$ to $h$ $w$ is contained in the $\varepsilon$-neigborhood of $H$. 
\endroster

(ii) If $A$ is a finite generating set for $G$ and $\delta >0$ is such that the geodesic triangles in $(\Gamma(G,A),d_A)$ are $\delta$-thin and $\lambda >0$ then any $1000\lambda\delta$-local $\lambda$-quasigeodesic is global $2\lambda$-quasigeodesic.

\endproclaim

\demo {Proof}
Part (i) of Proposition 3 easily follows from elementary properties of quasigeodesics and quasiconvex sets in hyperbolic spaces (see, for example [ABCFLMSS]) and from Theorem 3.4.4 of [ECHLPT] which asserts that any automatic structure consists of quasigeodesic words.
A careful proof is given in [Swa].
Part (ii) is established in [Gr, Lemma 7.2.B].
\enddemo
\smallskip
If any of the conditions (1)-(3) of part (i) of Proposition 3 is satisfied, we say that $H$ is {\it quasiconvex} in $G$.
\smallskip
Given a word hyperbolic group $G$ and its finite generating set $A$, the set $L$ of all $d_A$-geodesic words over $A$ is a regular language and is a part of an automatic structure for $G$. Suppose we already have this automatic structure at hand as well as the hyperbolicity constant $\delta$ of $(\Gamma (G,A), d_A)$ (that is geodesic triangles in $(\Gamma (G,A), d_A)$ are $\delta$-thin).
It follows from the discussion about quasiconvex subgroups of word hyperbolic groups in [BGSS] and [GS] that the knowledge of the quasiconvexity constant $\epsilon$ (and even the distortion parameter $C$)  for a quasiconvex subgroup $H$ of $G$ is sufficient for constructing the preimage $L_H$ of $H$ in $L$. 
Therefore the following statement is equivalent to Proposition 1 for word hyperbolic groups.
\proclaim {Proposition 4} Let $G$ be a word hyperbolic group given by a finite presentation $G=<x_1,..,x_n | r_1,..,r_m>$ and let $A=\{x_1^{\pm 1},..,x_n^{\pm 1}\}$. 

Then there is a uniform algorithm which,
given a finite set of words $v_1,..,v_t$ over $A$, will 
\roster
\item "(a)" eventually stop and produce the quasiconvexity constant $\epsilon$ and the distortion constant $C$ of the subgroup $H=gp(\overline {v_1},..\overline {v_t})$ of $G$ if $H$ is quasiconvex in $G$;
\item "(b)" run forever if $H$ is not quasiconvex in $G$.
\endroster
\endproclaim

\demo {Proof}

First we apply the Knuth-Bendix algorithm described in Section 6 of [ECHLPT] to produce the automaton $M$ recognizing the language $L=ShortLex(G,A)$  and the comparison automata $M_{=},M_{x_1},..,M_{x_n^{-1}}$ such that $$(L=L(M),A,M_{=},M_{x_1},..,M_{x_n^{-1}})$$ is an automatic structure with uniqueness for $G$. This allows us to calculate effectively in $G$.
Then we apply the algorithm of Papasoglu [P] or Olshansky [O] to produce the hyperbolicity constant $\delta$ of $\Gamma(G,A)$.
\smallskip
Put $V=(v_1,..,v_t,v_1^{-1},..,v_t^{-1})$. Let $K$ be the maximum of the lengths of words $v_1,..,v_t$.
\smallskip
{\bf Step 1} Start building the balls $N_1,N_2,..,N_i..$ of radius $i$ in the Cayley graph $\Gamma(H,V)$.
\smallskip
{\bf Step 2} For a current value of $i$ write down all $d_V$-geodesic words of length at most $2i$ and calculate minimal $\lambda>0$ such that the paths defined by these words in $\Gamma(G,A)$ are $\lambda$-quasigeodesics.
\smallskip
{\bf Step 3} Then for $j=i+1,..,1000Ki\delta\lambda$ check if all $d_V$-geodesic words of length at most $2j$ are still $\lambda$-quasigeodesics in $\Gamma(G,A)$.
If for some of these $j$ it is not so, we go back to Step 1 and increase $i$ by 1.
If for each of these $j$ we succeed, then $H$ is quasiconvex in $G$.
Indeed, in this case any $d_V$-geodesic word defines a $1000\delta\lambda$-local $\lambda$-quasigeodesic in $\Gamma(G,A)$ and therefore by part (ii) of Proposition 3 this path is $2\lambda$-quasigeodesic in $\Gamma(G,A)$.
Part (i) of Proposition 3 implies that $H$ is quasiconvex in $G$ with the distortion constant $2\lambda$.
\smallskip
On the other hand, if $H$ is quasiconvex in $G$ then for some $\lambda>0$ all $d_V$-geodesics define $\lambda$-quasigeodesics in $\Gamma(G,A)$ and so our algorithm will eventually disclose this fact and stop.
\smallskip
The promised quasiconvexity constant $\epsilon$ for $H$ can be taken to be $1000\delta(1+log_2 C)$ where $C$ is the distortion constant (see [Gr, Proposition 7.2.A]).

\enddemo

{\bf Example} Let $k\ge 2$ and $$G=<a_1,b_1,..,a_k,b_k | \prod_{i=1}^k [a_i,b_i]>$$ and $v_1,..,v_t$ be some words over $a_1,..,b_k$ defining the subgroup $H=gp(\overline {v_1},..,\overline {v_t})$ of $G$.

 It is well known (see, for example [Swa]) that the surface group $G$ is word hyperbolic and $H$ is quasiconvex in $G$. Thus the generalized word problem for $G$ with respect to $H$ is solvable.
We are not so much interested in the reasons for the quasiconvexity and even in finding a presentation for $H$ but we still can apply algorithms described in Proposition 1 or Proposition 4 to solve the generalized word problem.
Moreover, we can practically solve the {\it generation problem} for $G$, that is, decide if a given finite collection of elements of $G$ generate $G$ or not. To do this we first construct an automatic structure with uniqueness $L=ShortLex(G)$ for $G$  and then find the preimage $L_H$ in $L$ of the subgroup $H$, generated by this collection of elements.
It remains to check if $L-L_H$ is empty or not. 
\smallskip
{\bf Remark} Notice that for the class of automatic groups solvable generalized word problem does not imply quasiconvexity. Let $G$ be as in the example above and $\phi$ be an automorphism of $G$ coming from a pseudo-anosov homeomorphism of the surface whose fundamental group is $G$.
Take $G_1$ to be the HNN-extension of $G$ along $\phi$ that is
$$G_1=<a_1,b_1,..,a_k,b_k,t| \prod_{i=1}^k [a_i,b_i], t^{-1}a_1t=\phi(a_1),..,t^{-1}b_kt=\phi(b_k)>$$
Then (see [Th]) $G_1$ is word hyperbolic and there is a short exact sequence $$1\rightarrow G\rightarrow G_1\rightarrow {\Bbb Z}\rightarrow 1.$$ 
Since $G$ and ${\Bbb Z}$ are infinite, by Proposition 3.9 of [ABCFLMSS] $G$ is not quasiconvex in $G_1$. This also implies that if $L$ is an automatic language for $G_1$ then $G$ is not $L$-rational (see [Swa] for details). On the other hand the generalized word problem for $G_1$ with respect to $G$ is solvable: if $v$ is a word in the generators of $G_1$ then it represents an element of $G$ if and only if the letter $t$ occurs in $v$ with the exponent sum zero. 
\smallskip
Notice also that there is an asynchronously automatic structure on $G_1$ for which $G$ is rational (see [Sha] for details). At the moment we do not know whether there are finitely generated subgroups of asynchronously automatic groups with solvable generalized word problem which are not rational for all asynchronously automatic structures on the ambient groups. 
\smallskip
We also do not know if there is an automatic group $G$ and a finitely presentable subgroup $H$ such that $G$ has unsolvable generalized word problem with respect to $H$.
E.Rips [Ri] constructed an example of a finitely generated (but not finitely presented) subgroup $H$ of a word hyperbolic group $G$ with unsolvable generalized word problem.

\enddocument